% Paper:        " On real analytic Banach manifolds "
%
% Date:         October 10, 2007
% Authors:       (Mr.) Imre PATYI
%               Department of Mathematics and Statistics
%               Georgia State University
%               Atlanta, GA 30303-3083
%               USA
% Email:        ipatyi@gsu.edu
% Tel. 404-413-6450
% Fax. 404-413-6403
%               Scott Bradford Simon
%               Department of Mathematics
%               SUNY Stony Brook
%               Stony Brook, NY 11794-3651
%               USA
% Email:        sbsimon@math.sunysb.edu
% Tel. 631-632-4506
% Fax. 631-632-7631
% Remarks on the file:

% 0.  This is an AmS-TeX file, not AmS-LaTeX.
% 1.  There are many abbreviations and a few macros.
% 2.  If any difficulty arises with the file or the
%    content of this manuscript, please email the author
%    at "ipatyi@gsu.edu".
% 3.  This file contains 1600 lines.

% PAPER BEGIN
% ~~~~~~~~~~~~~~~~~~~~~~~~~~~~~~~~~~~~~~~~~~~~~~~~~~~~~~~~~~~~~~~~~~~~~~~~~~~~~
% ~~~~~~~~~~~~~~~~~~~~~~~~~~~~~~~~~~~~~~~~~~~~~~~~~~~~~~~~~~~~~~~~~~~~~~~~~~~~~
% ~~~~~~~~~~~~~~~~~~~~~~~~~~~~~~~~~~~~~~~~~~~~~~~~~~~~~~~~~~~~~~~~~~~~~~~~~~~~~
% ~~~~~~~~~~~~~~~~~~~~~~~~~~~~~~~~~~~~~~~~~~~~~~~~~~~~~~~~~~~~~~~~~~~~~~~~~~~~~

%\input /math/faculty3/ipatyi/tex/amstex.tex

\input amstex.tex

\magnification=\magstep1
\hsize=5.5truein
\vsize=9truein
\hoffset=0.5truein
\parindent=10pt
\newdimen\nagykoz
\newdimen\kiskoz
\nagykoz=7pt
\kiskoz=2pt
\parskip=\nagykoz
\baselineskip=12.7pt

%\raggedbottom

\loadeufm \loadmsam \loadmsbm

\font\vastag=cmssbx10
\font\drot=cmssdc10
\font\vekony=cmss10
\font\vekonydolt=cmssi10
\font\cimbetu=cmssbx10 scaled \magstep1
\font\szerzobetu=cmss10

\font\scVIII=cmcsc8
\font\rmVIII=cmr8
\font\itVIII=cmti8
\font\bfVIII=cmbx8
\font\ttVIII=cmtt8

\def\cim#1{{\centerline{\cimbetu#1}}}
\def\szerzo#1{{\vskip0.3truein\centerline{\szerzobetu#1}}}
\def\alcim#1{{\medskip\centerline{\vastag#1}}}
\def\tetel#1#2{{{\drot#1}{\it\szukebb~#2\tagabb}}}
\long\def\biz#1#2{{{\vekony#1} #2}}
\def\kiemel#1{{\vekonydolt#1\/}}
\long\def\absztrakt#1#2{{\vskip0.4truein{\vekony#1} #2\vskip0.5truein}}
\def\szukebb{\parskip=\kiskoz}
\def\tagabb{\parskip=\nagykoz}
\def\vonal{{\vrule height 0.2pt depth 0.2pt width 0.5truein}}

\def\CC{{\Bbb C}}
\def\PP{{\Bbb P}}

\def\bmfd{{Banach manifold}}
\def\bsmfd{{Banach submanifold}}
\def\hmfd{{Hilbert manifold}}
\def\hsmfd{{Hilbert submanifold}}
\def\balg{{Banach algebra}}
\def\blgp{{Banach Lie group}}
\def\bla{{Banach Lie algebra}}
\def\bvbdl{{Banach vector bundle}}
\def\hvbdl{{Hilbert vector bundle}}

\def\Ker{\hbox{\rm Ker}}

\def\im{\hbox{\rm Im}}

\def\cC{{\Cal C}}

\def\Oa{{\Omega}}
\def\oa{{\omega}}
\def\UU{{\frak U}}

\def\GL{\text{\rm GL}}

\def\Hom{\text{\rm Hom}}
\def\End{\text{\rm End}}

\def\Aut{\text{\rm Aut}}

\def\vbdl{{vector bundle}}
\def\hspc{{Hilbert space}}

\def\cts{{continuous}}
\def\bdd{{bounded}}

\def\idml{{infinite dimensional}}

\def\idms{{infinite dimensions}}

\def\fdml{{finite dimensional}}

\def\pscx{{pseudoconvex}}

\def\ubs{{unconditional basis}}

\def\st{{such that}}
\def\wrt{{with respect to}}

\def\delbar{{\bold{\bar\partial}}}

\def\<{{\langle}}
\def\>{{\rangle}}
\def\NN {{\Bbb N}}
\def\RR {{\Bbb R}}

\def\AA{{\Cal A}}

\def\OO {{\Cal O}}

\def\fii{\varphi}
\def\fn{func\-tion}
\def\fns{func\-tions}
\def\holo{hol\-o\-mor\-phic}
\def\ra{real analytic}

\def\mfd{manifold}
\def\smfd{submanifold}
\def\cpx{complex}
\def\cpt{compact}

\def\nbd{neighbor\-hood}

\def\bspc{Banach space}

\def\La{\Lambda}

\def\da{{\delta}}
\def\Da{{\Delta}}

\def\da{{\delta}}

\def\Prop{Proposition}
\def\p#1{{\Prop~#1}}

\def\Th{{Theorem}}
\def\th{{theorem}}

\def\t#1{{\Th~#1}}

% today.tex: Macro to print today's date
%--
% Author: John W. Shipman, NM Tech Computer Center,
%   Socorro, NM 87801; john at nmt.edu
%--
% EXPORTED FUNCTIONS:
%   \today: Outputs today's date as ``yyyy-mm-dd''
%   \now: Outputs the current time as ``hh:dd''
%   \timestamp: \today, plus one space, plus \now
%--
\newcount\minute    % Current minute within the hour
\newcount\hour      % Current hour (24-hour type)
\newcount\hourMins  % Temporary for taking hour modulo 60
%
% - - -   \ n o w   - - -
%
\def\now%
{% Displays today's time as ``hh:mm''
%  The \time macro gives the minutes since midnight.  Compute
%  the whole hours by dividing this by 60, then find the
%  minute by effectively taking the minutes modulo 60.
%
  \minute=\time    % Number of minutes since midnight
  \hour=\time \divide \hour by 60 % Get hours
  \hourMins=\hour \multiply\hourMins by 60
  \advance\minute by -\hourMins % Hours modulo 60
  \zeroPadTwo{\the\hour}:\zeroPadTwo{\the\minute}%
}% --- \now ---
%
% - - -   \ t i m e s t a m p   - - -
%
\def\timestamp%
{% Displays ``yyyy-mm-dd hh:mm''
  \today\ \now
}% --- \timestamp ---
%
% - - -   \ t o d a y   - - -
%
\def\today%
{% Displays today's date and time as ``yyyy-mm-dd hh:mm''
  \the\year-\zeroPadTwo{\the\month}-\zeroPadTwo{\the\day}%
}% --- \today ---
%
% - - -   z e r o P a d T w o   - - -
%
\def\zeroPadTwo#1%
{% Left zero pad of the argument to 2 digits.  The argument
%  should be a number between 1 and 99.  This macro outputs
%  a `0' if the argument is less than ten, then it outputs
%  the argument.
%
  \ifnum #1<10 0\fi    % Conditionally output a zero
  #1%    Then output the argument
}% --- \zeroPadTwo ---

%{\ttVIII\jobname, \timestamp}
{\phantom.}
\vskip0.5truein
\cim{ON REAL ANALYTIC BANACH MANIFOLDS}
%\vskip0.2truein
%\cim{BUNDLES OVER BANACH SPACES}
\szerzo{Imre Patyi and Scott Bradford Simon\plainfootnote{${}^*$}{\rmVIII 
Both authors were supported in part by NSF grants 
DMS 0600059, DMS 0203072, and a
Research Initiation Grant from Georgia State University.}}
\absztrakt{ABSTRACT.}{Let $X$ be a real Banach space with an unconditional basis
        (e.g., $X=\ell_2$ Hilbert space), $\Omega\subset X$ open,
        $M\subset\Omega$ a closed split real analytic Banach submanifold of $\Omega$,
        $E\to M$ a real analytic Banach vector bundle, and ${\Cal A}^E\to M$ the sheaf of 
	germs of real analytic sections of $E\to M$.
         We show that the sheaf cohomology groups $H^q(M,{\Cal A}^E)$ vanish for all $q\ge1$,
        and there is a real analytic retraction $r:U\to M$ from an open set $U$ with
        $M\subset U\subset\Omega$ such that $r(x)=x$ for all $x\in M$.
         Some applications are also given, e.g., we show that any infinite dimensional
	real analytic Hilbert submanifold of separable affine or projective Hilbert space
	is real analytically parallelizable.



%Running head: Real analytic Banach manifolds, Patyi and Simon.

}

% Section numbering
\def\sA{{1}}
\def\sB{{2}}
\def\sC{{3}}
\def\sD{{4}}
\def\sE{{5}}
\def\sF{{6}}
\def\sG{{7}}
\comment

\endcomment

% Proclaim numbering (theorems) \tAA etc.
% Display numbering (equations) \eAA etc.

\def\tAA{{\sA.1}}

\def\tBA{{\sB.1}}
\def\tBB{{\sB.2}}

\def\tCA{{\sC.1}}
\def\tCB{{\sC.2}}
\def\tCC{{\sC.3}}

\def\tEA{{\sE.1}}
\def\tEB{{\sE.2}}
\def\tEC{{\sE.3}}

\def\tFA{{\sF.1}}
\def\tFB{{\sF.2}}
\def\tFC{{\sF.3}}

\def\tGA{{\sG.1}}
\def\tGB{{\sG.2}}
\def\tGC{{\sG.3}}

\comment

\endcomment

% Reference numbering
\def\rC{C}
\def\rD{D}
\def\rDPV{DPV}
\def\rE{E}
\def\rF{F}
\def\rG{G}
\def\rJE{J5}
\def\rH{H}
\def\rK{K}
\def\rLA{L1}

\def\rLC{L2}
\def\rLD{L3}
\def\rLE{L4}
\def\rLP{LP}
\def\rLg{Lg}
\def\rM{M}
\def\rPA{P1}
\def\rPB{P2}
\def\rPC{P3}
\def\rS{S}
\def\rWB{WB}
\def\rZn{Zn}
\def\rZg{Zg}

\alcim{\sA. INTRODUCTION.}

         Many of the classical \mfd{}s that occur in mathematics are \ra\ \mfd{}s,
        not just $C^\infty$-smooth \mfd{}s.
         A natural question over a \ra\ \mfd\ is the existence and abundance of
        global \ra\ \fns\ and \ra\ sections of \ra\ \vbdl{}s.
         It is natural to expect and quite well-known that \ra\ \fns\ and
        \ra\ \mfd{}s are intimately related to \holo\ \fns\ and \cpx\ \mfd{}s.

         In the 1950s a powerful theory of \cpx\ \mfd{}s, called Stein theory,
        was developed, which to a large degree answers the question of existence
        and abundance of global \holo\ \fns\ and \holo\ sections of \holo\ \vbdl{}s
        over Stein manifolds, i.e., over closed \smfd{}s of \cpx\ Euclidean spaces.
         Among the first applications of this Stein theory was a study of the
        question of existence and abundance of global \ra\ \fns\ and \ra\ sections
        of \ra\ \vbdl{}s over \ra\ \mfd{}s.
         This was done in two phases.
         First for \ra\ \smfd{}s of Euclidean spaces, then for arbitrary abstract
        \ra\ \mfd{}s.
         There seem to be two classical routes to the problem of existence and
        abundance over a \ra\ \mfd.
         The first one is based on complexification and application of results from
        the theory of Stein \mfd{}s.
         The second one relies on the Malgrange approximation theorem, and draws
        on techniques from analysis that cannot currently be matched in \idms.

         Here we follow the method of complexification and apply some results
        from an emerging theory of \cpx\ \bsmfd{}s of \bspc{}s mainly developed
        by Lempert and some of his students.
         Our main point is \t\tAA\ below.

\tetel{\t\tAA.}{Let $X$ be a real \bspc\ with an \ubs\
        (e.g., $X=\ell_2$ Hilbert space), $\Oa\subset X$ open,
        $M\subset\Oa$ a closed split \ra\ \bsmfd\ of\/ $\Oa$,
        $E\to M$ a \ra\ \bvbdl\ with a \bspc\ $Z$ for fiber type, 
	and $\AA^E\to M$ the sheaf of germs of \ra\
        sections of $E\to M$.
         Then the following hold.
\vskip0pt
	{\rm (a)} There is a \ra\ map $r:\oa\to M$ from an open set $\oa$ with
        $M\subset\oa\subset\Oa$ \st\ $r(x)=x$ for all $x\in M$; i.e., $M$ is a
        \ra\ \nbd\ retract in\/ $\Oa$.
\vskip0pt
        {\rm (b)} The sheaf cohomology groups $H^q(M,\AA^E)$ vanish for
        all $q\ge1$.
\vskip0pt
        {\rm (c)} If $E\to M$ is topologically trivial, then $E\to M$ is
        \ra{}ally trivial.
\vskip0pt
        {\rm (d)} If $Z$ is isomorphic to the Hilbert space $\ell_2$ or
        $M$ is contractible, then $E\to M$ is \ra{}ally trivial.
\vskip0pt
        {\rm (e)} For $1\le p<\infty$ let
        $Z_1=\ell_p(Z)=\{z=(z_n):z_n\in Z,\ \|z\|=(\sum_{n=1}^\infty\|z_n\|^p)^{1/p}<\infty\}$.
        The \bvbdl\ $E\oplus(M\times Z_1)\to M$ is \ra{}ally isomorphic to the
        trivial bundle $M\times Z_1$.
}

         Some applications of \t\tAA\ are also given in \S\S\,\sE-\sG.

\alcim{\sB. BACKGROUND.}

         In this section we collect some definitions and theorems that are useful for
        this paper.  
	 Some good sources of information on \cpx\ analysis on \bspc{}s are 
	[\rD, \rM, \rLA].

	 Put $B_X(x_0,r)=\{x\in X:\|x-x_0\|<r\}$ for the open ball with center $x_0\in X$
	and radius $0<r\le\infty$ in any real or \cpx\ \bspc\ $(X,\|\cdot\|)$, and
	let $B_X(r)=B_X(0,r)$.

         A \kiemel{\cpx\ \bmfd} $M$ modelled on a \cpx\ \bspc\ $X$ is a paracompact
        Hausdorff space $M$ with an atlas of bi\holo{}ally related charts onto
        open subsets of $X$.
         A subset $N\subset M$ is called a \kiemel{closed \cpx\ \bsmfd} of $M$
        if $N$ is a closed subset of $M$ and for each point $x_0\in N$ there are
        an open \nbd\ $U$ of $x_0$ in $M$ and a bi\holo\ map $\fii:U\to B_X(1)$ onto
        the unit ball $B_X(1)$ of $X$ that maps the pair $(U,U\cap N)$ to
        a pair $(B_X(1),B_X(1)\cap Y)$ for a closed \cpx\ linear subspace $Y$ of $X$.
         The \smfd\ $N$ is called a \kiemel{split} or \kiemel{direct} \bsmfd\ of $M$
        if at each point $x_0\in N$ the corresponding subspace $Y$ has a direct complement
        in $X$.
	 A \kiemel{maximally real submanifold} $N'$ of a complex manifold $M$ is a
	real submanifold such that the \holo\ tangent space $T^{1,0}_x M$ is 
	the complexification of $T_x N'$ for every $x \in N'$.  
	 Equivalently, there is an atlas $(U_\alpha)$ on $M$ mapping 
	each $U_\alpha \cap N'$ to a real subspace of $X$ whose complexification is $X$.

	 We call a map $f:U\to V$ from an open subset of a real \bspc\ onto
	an open subset of a real \bspc\ \kiemel{bi-\ra} or a \kiemel{\ra\
	diffeomorphism with \ra\ inverse} if both $f$ and its inverse
	$f^{-1}:V\to U$ are \ra.
	 A \kiemel{\ra\ \bmfd} $M$ modelled on a real \bspc\ $X$ is
	 a paracompact
	Hausdorff space $M$ with an atlas of bi-\ra{}ally related charts
	onto open subsets of $X$.
	 A subset $N\subset M$ is called a \kiemel{closed \ra\ \bsmfd}
	of $M$ if $N$ is a closed subset of $M$ and for each point $x_0\in
	N$ there are an open \nbd\ $U$ of $x_0$ in $M$ and a bi-\ra\
	map $\fii:U\to B_X(1)$ onto the unit ball $B_X(1)$ of $X$ that
	maps the pair $(U,U\cap N)$ to a pair $(B_X(1),B_X(1)\cap Y)$
	for a closed real linear subspace $Y$ of $X$.  The \smfd\ $N$
	is called a \kiemel{split} or \kiemel{direct} \bsmfd\ of $M$
	if at each point $x_0\in N$ the corresponding subspace $Y$ has
	a direct complement in $X$.

         We call a subset $S$ of a \cpx\ \bspc\ $X$ endowed with a complex conjugation
	$X\to X$, $x\mapsto\bar x$,
        \kiemel{symmetric \wrt\ complex conjugation} if $\bar x\in S$ for all $x\in S$.
	For \cpx\ \bspc{}s $Z,Z'$ denote by $\Hom(Z,Z')$ the \bspc\ of all
	\bdd\ \cpx\ linear operators $A:Z\to Z'$ endowed with the operator norm.
	Write $\End(Z)$ for $\Hom(Z,Z)$ and $\GL(Z)$ for the group of
	invertible elements of the \balg\ $\End(Z)$; $\GL(Z)$ is a \blgp\
	with a \bla\ $\End(Z)$.
	If $Z$ has a \cpx\ conjugation, then $\End(Z)$ can also be endowed
	with the induced \cpx\ conjugation defined by $A\mapsto\bar A$,
	$\overline{A\overline{z}}=\bar Az$.
	Then $\overline{AB}=\bar A\bar B$ for $A,B\in\End(Z)$, and
	$A\in\End(Z)$ is invertible if and only if $\bar A$ is.

	 We denote the spaces or sheaves of all \holo\ \fns\ or sections
	by $\OO$, \ra\ ones by $\AA$, and \cts\ ones by $\cC$.

	 The following \th\ about \cpx\ \bmfd{}s will come in handy.

\tetel{\t\tBA.}{Let $X$ be a \cpx\ \bspc\ with a \cpx\ conjugation,
	$\Oa\subset X$ \pscx\ open, and $M\subset\Oa$ a closed split \cpx\
	\bsmfd\ of\/ $\Oa$.
	 Suppose that $X$ has an \ubs, and $M$ is symmetric \wrt\ \cpx\ conjugation.
	 Then the following hold.
\vskip0pt
	{\rm(a)} There are a \pscx\ open set $\oa\subset X$ and a \holo\
	map $r:\oa\to M$ \st\ $M\subset\oa\subset\Oa$, $\oa$ is symmetric \wrt\
	\cpx\ conjugation, $\overline{r(\bar x)}=r(x)$ for all $x\in\oa$, and
	$r(x)=x$ for all $x\in M$; i.e., $M$ admits a real-type \holo\ retraction $r$
	from one of its \pscx\ open \nbd{}s $\oa$.
\vskip0pt
	{\rm(b)} Let $Z$ be a \cpx\ \bspc\ with a \cpx\ conjugation,
	$G=\GL(Z)$ endowed with the \cpx\ conjugation induced from that of $Z$,
	$\UU$ an open covering of\/ $\Oa$, and $(g_{UV})\in Z^1(\UU,\OO^G)$ a \holo\
	cocycle.
	 If\/ $\Oa$ and all members $U$ of\/ $\UU$ are symmetric \wrt\ \cpx\ conjugation,
	there is a \cts\ cochain $(c_U)\in C^0(\UU,\cC^G)$ with
	$c_U(x)^{-1}g_{UV}(x)c_V(x)=1$ for all $x\in U\cap V$, $U,V\in\UU$, and
	$g_{UV}$ and $c_U$ are of real-type in the sense that
	$\overline{g_{UV}(\bar x)}=g_{UV}(x)$ and
	$\overline{c_{U}(\bar x)}=c_U(x)$, then
	there is a \holo\ cochain $(d_U)\in C^0(\UU,\OO^G)$ \st\
	$d_U(x)^{-1}g_{UV}(x)d_V(x)=1$ for all $x\in U\cap V$, 
	and $\overline{d_U(\bar x)}=d_U(x)$ for all $x\in U$, $U,V\in\UU$.
}

\biz{Proof.}{The analogs of (a) and (b) without special attention to \cpx
	conjugation were proved in general in [\rLP] for (a) and [\rPC] for (b).
	 The (otherwise fairly long and involved) proofs there admit
	straightforward modifications to take into account \cpx\ conjugations,
	thereby proving \t\tBA.
}

	 We shall also make use of the \th\ below,
	whose formulation is eased by the following definition.
	 Let $M$ be a \cpx\ \bmfd, $E\to M$ a \holo\ \bvbdl,
	$E^{p,q}\to M$ the \ra\ \bvbdl\ of $(p,q)$-forms with values in
	$E$ over $M$ for $p,q\ge0$.
	 Note that $E^{p,0}\to M$ is a \holo\ \bvbdl.
	 Consider the \ra\ Dolbeault groups $H^{p,q}_{\delbar,\oa}(M,E)$
	of $M$ with coefficients in $E$ for $p,q\ge0$ defined as
$$
	H^{p,q}_{\delbar,\oa}(M,E)=\frac{\Ker(\delbar:C^\oa(M,E^{p,q})\to C^\oa(M,E^{p,q+1}))}
		{\im(\delbar:C^\oa(M,E^{p,q-1})\to C^\oa(M,E^{p,q}))},
 $$
 	where $\im=0$ for $q=0$.
	 We say that 
	\kiemel{the \ra\ Dolbeault isomorphism holds over $M$}
	if the \ra\ Dolbeault groups $H^{p,q}_{\delbar,\oa}(M,E)$ are canonically
	isomorphic to the sheaf cohomology groups $H^q(M,\OO^{E^{p,0}})$
	in the usual way for $p,q\ge0$.

\tetel{\t\tBB.}{{\rm(S.\,B.\,Simon, [\rS])}
	{\rm(a)} Let $X$ be a real \bspc\ with an \ubs, and\/ $\Oa\subset X$ open.
	 If $U\subset X'$ is any open subset of the complexification $X'$ of $X$
	with\/ $\Oa\subset U$, then there is a \pscx\ open\/ $\Oa'\subset X'$ with\/
	$\Oa=\Oa'\cap X$ and\/ $\Oa'\subset U$; i.e., $\Oa\subset X$ has a \nbd\
	basis consisting of \pscx\ open subsets\/ $\Oa'$ of $X'$.
\vskip0pt
	{\rm(b)} Let $X,\Oa$ be as in\/ {\rm(a)}, and $E\to\Oa$ a \ra\ \bvbdl.
	 Then the sheaf cohomology groups $H^q(\Oa,\AA^E)$ of \ra\ sections 
	vanish for all $q\ge1$.
\vskip0pt
	{\rm(c)} If $X$ is a \cpx\ \bspc\ with an \ubs,
	then the \ra\ Dolbeault isomorphism holds over any open\/ $\Oa\subset X$.
}

\alcim{\sC. COMPLEXIFICATION.}

         In this section we see how to complexify a closed \ra\ \bsmfd\ of an 
	open subset of a Banach space, as well as an abstract \ra\ \bmfd.

\tetel{\p\tCA.}{Let $M_i$ be real analytic, maximally real Banach
	submanifolds of complex Banach manifolds $M_i'$ for $i=1,2.$
	 Suppose $f:M_1\to M_2$ is bi-\ra\ onto\/ $M_2$.
	 Then $f$ extends to a bi\holo\ \fn\ from an open subset of
	$M'_1$ onto an open subset of $M'_2.$
}

\biz{Proof.}{There is a locally finite atlas $(U_{i,\alpha}, \phi_{i,\alpha})$ for $M'_i$,
	where for each $i,\alpha$ there is a real Banach space $X_{i,\alpha}$ 
	with complexification $X_{i, \alpha}'$ such
	that $\phi_{i,\alpha}:U_{i,\alpha} \rightarrow X'_{i,\alpha}$ is a 
	biholomorphism onto an open subset of $X'_{i, \alpha}$ and
	$\phi_{i,\alpha}(U_{i,\alpha} \cap M_i) \subset X_{i,\alpha}$.
	
	 Let $g=f^{-1}$.  
	 After shrinking each $U_{i,\alpha},$ we may
	assume that $f$ extends holomorphically to a function $f_{1,\alpha}$ 
	on each $U_{1,\alpha}$, and likewise for $g$ to $g_{2,\alpha}$ on each $U_{2,\alpha}$.

	 We will refine these covers so that the intersections are connected:  for each
	$x \in M_i \cap U_{i,\alpha},$ choose $U_i(x) \subset
	U_{i,\alpha}$ open such that $U_i (x) \subset U_{i,\beta}$
	whenever $x \in U_{i,\beta},$ and $U_i(x) \cap U_{i,\gamma} \ne
	\emptyset$ if and only if $x \in \overline{U}_{i,\gamma}$.  
	 Now take 
$$
	V_{i,\alpha}=\bigcup_{x \in U_{i,\alpha} \cap M_i} U_i (x).
 $$
	 This forms a refinement of $(U_{i,\alpha})$ whose
	intersections are connected, and each intersection intersects
	$M_i$ as well.  Note that this cover may no longer cover all of $M_i'$.

	 Let $V_i=\bigcup_\alpha V_{i,\alpha}$.  
	 By the uniqueness of holomorphic extension of real analytic functions 
	(see [\rS, Proposition 6.1]) the functions $f_{1,\alpha}$, respectively
	$g_{2,\alpha}$ agree on the overlaps, and thus define functions 
	$f_1$, respectively $g_2$ on $V_1$, respectively $V_2$.  
	 Let $W_1=V_1\cap g_2(V_2),$ and let $W_2=f_1(W_1)$.
	 Observe that $f_1(W_1) \subset V_2$ and 
	so $g_2 \circ f_1$ and $f_1 \circ g_2$ are the respective identity functions 
	on $W_1$ and $W_2$. 
	 Thus $f_1$ is a biholomorphism on $W_1$.
	 The proof of \p\tCA\ is complete.
}

	 The next theorem states that it is possible to complexify 
	abstract real analytic Banach manifolds.

\tetel{\t\tCB.}{Let $M$ be a real analytic Banach manifold.  
	Then there is a complex manifold $M'$ such that $M$ is a
	maximally real submanifold of $M'$.
	 Furthermore, if $(U_\alpha, \phi_\alpha)$ is an atlas of $M$,
	there is an atlas $(V'_a, \psi_a)$ of $M'$ whose
	restriction to $M$ is a refinement of $(U_\alpha, \phi_\alpha)$.
}

\biz{Proof.}{Let $(U_\alpha, \phi_\alpha)$ be a locally finite atlas for
	$M$, and let $(V_a)$ be a locally finite strict refinement, i.e., if
	$I$ is the index set for $(U_\alpha)$ and $J$ is the index set for
	$(V_a)$, then there is a map $\sigma: J \rightarrow I$ such that
	$\overline{V}_a \subset U_{\sigma a}$. 
	 Let $X_\alpha$ be the target space
	for $\phi_\alpha$.
	 For each $\alpha, \beta \in I$, 
	the map $\phi_{\beta \alpha}=\phi_\beta \circ \phi_\alpha^{-1}$
	is bi-real analytic where defined. 

	 Define $U_{\beta \alpha}=\phi_\alpha(U_\alpha \cap U_\beta)$ and likewise
	$V_{ba}=\phi_{\sigma a}(V_a \cap V_b)$. 
	 Let $\Re$ be the projection taking a complex vector to its real part,
        and let $X_\alpha'$ be the complexification of $X_\alpha$.
	 There is a neighborhood $U'_{\beta \alpha} \subset X_\alpha'$ of $U_{\beta \alpha}$ 
	such that $\phi_{\beta \alpha}$ extends biholomorphically to a function
	$\psi_{\beta \alpha}$ on $U'_{\beta \alpha}$ (see \p\tCA). 
	 Without loss of generality, we may assume that for every $\alpha, \beta \in I$,
	every component of 
	$U'_{\beta \alpha} \cap \psi_{\alpha \beta}(U'_{\alpha \beta})$ intersects
	$U_{\alpha \beta}$.  
	 Observe that $\psi_{\alpha \beta} \circ \psi_{\beta \alpha}$ and
	$\psi_{\alpha \beta} \circ \psi_{\beta \gamma} \circ \psi_{\gamma \alpha}$ 
	are both the identity where defined, by uniqueness of holomorphic extensions.  
	 Define
	$U'_{\gamma \beta \alpha}=U'_{\gamma \alpha}\cap
	\psi_{\alpha \beta}(U'_{\gamma \beta}\cap U'_{\alpha \beta})$ 
	and 
	$V_{cba}=\phi_{\sigma a} (V_a \cap V_b \cap V_c)$.
	  For each $p \in V_a$, define $V'_a (p) \subset X'_{\sigma a}$ to be a
	connected neighborhood of $\phi_{\sigma a} (p)$ such that
$$
\aligned
\text{(1) }& \Re V'_a (p) \subset V'_a(p) \cap \phi_{\sigma a} (V_a), \\
\text{(2) }& V'_a(p) \cap \psi_a (V'_b) \ne \emptyset 
	\text{ if and only if}\\
	& p \in \overline{V_b}
	\text{ (of which there are only finitely many),} \\
\text{(3) }& \text{if } p \in \overline{V}_a \cap \overline{V}_b, \text{ then } V'_a (p) \subset U_{\sigma b \sigma a}, \text{and} \\
\text{(4) }& V'_a(p) \cap U'_{\gamma \beta \sigma a} \subset U_{\gamma \sigma a}'
	\text{ for all } \beta, \gamma \in I.
\endaligned
 $$
	 Note that since 
	$\phi_{\sigma a} (V_a) \cap U'_{\gamma \beta \sigma a} \subset U_{\gamma \sigma a}'$ 
	and because of the local finiteness of the cover $(U_\alpha')$,
	(4) is possible to satisfy.
	 Let
$$
	 V'_a=\bigcup_{p \in V_a} V'_a(p).
 $$
	 By (1), the collection of sets $(V_a')$ is locally finite.
	 We say that $x\sim y$ if $x \in \overline{V}'_a$ for some $a$ and
	$y \in \overline{V}'_b$ for some $b$ and $\psi_{\sigma b \sigma
	a}(x)=y$, where $\psi_{\alpha \alpha}$ is the identity. 

	 This is an equivalence relation: 
	since $\psi_{\alpha \beta}= \psi_{\beta \alpha}^{-1}$, 
	this relation is symmetric.  
	 Clearly, it is reflexive. 
	 For transitivity, suppose $x\sim y$ and $y\sim z$.  
	 Then there are $a, b$, and $c\in J$ such that
	$\psi_{\sigma b\sigma a}(x)=y$ and $\psi_{\sigma c\sigma b}(y)=z$.  
	 We wish to show that $x$ is in the domain of $\psi_{\sigma c \sigma a}$ 
	and $z$ is in its image.
	 Transitivity will follow from the fact that
	$\psi_{\sigma c \sigma b} \circ \psi_{\sigma b \sigma a}= \psi_{\sigma c \sigma a}$ 
	where they are all defined. 
	 This last is true on the real subspace, and therefore on the complexification.  
	 Without loss of generality, $x \in V_a' (p)$ for some $p \in V_a$.  
	 Since $\psi_{\sigma b\sigma a}(x)=y$, we have that $y \in U_{\sigma a \sigma b}'$, 
	and since $\psi_{\sigma c\sigma b}(y)=z$,
	we have that $y \in U_{\sigma c \sigma b}'$.
	 Therefore, 
	$x \in \psi_{\sigma a \sigma b}(U_{\sigma a \sigma b}' \cap U_{\sigma c \sigma b}')$.
	 By (4), we conclude $x \in U_{\sigma c \sigma a}'$.  
	 To show that $z \in U_{\sigma a \sigma c}'$ we interchange the roles
	of $x$ and $z$ and the roles of $a$ and $c$.  
	 Thus, transitivity is proved.

	 Observe that if $x \in V_a, y \in V_b,$ then $x \sim y$ if and only
	if $\phi_{\sigma a}^{-1}(x)=\phi_{\sigma b}^{-1}(y).$

	 Thus, we can identify $M$ with its image under $\phi$ modulo $\sim,$ and we can take
$$
	M'=\Big(\bigcup_{a \in J} V'_a \Big)/\sim.
 $$
	 We must check that $M'$ is Hausdorff.  
	 Let $p,q \in M',$ and suppose every neighborhood of $p$ intersects 
	every neighborhood of $q$.
	 Let $x \in V'_a$ be a representative of the equivalence class
	defined by $p$, and $y \in V'_b$ a representative of the
	equivalence class defined by $q$.  
	 Choose $p' \in V_a$ and $q' \in V_b$ such that $x \in V'_a(p')$ 
	and $y \in V'_b(q')$.  
	 Then $V'_a(p') \cap \psi_{\sigma a \sigma b}(V'_b(q'))$ is nonempty,
	and therefore by (2) and (3), $V'_a(p') \subset U'_{\sigma b \sigma a}$ 
	and $V'_b(q') \subset U'_{\sigma a \sigma b}$. 

	 This tells us that both $x$ and $\psi_{\sigma a \sigma b}(y)$ are
	in $U_{\sigma a}$, and every neighborhood from each intersects the other.  
	 By the fact that the complexification $X_{\sigma a}'$ of $X_\alpha$ is Hausdorff,
	we can see that $x = \psi_{\sigma a \sigma b}(y)$.
	 It is easy to see that $T_x^{1,0} M'$ is the complexification 
	of $T_x M$ for each $x \in M$, so $M \subset M'$ is maximally real.
	 The proof of \t\tCB\ is complete.
}

\tetel{\t\tCC.}{Let $M$ be a real analytic Banach manifold, 
	$N\subset M$ a closed real analytic Banach submanifold.
         Then there are a complex Banach manifold $M'$ and $N'\subset M'$ a 
	 closed \cpx\ \bsmfd\ of $M'$ \st\/ $M \subset M'$ and $N \subset N'$ are maximally
         real.  
	  In particular, if $X$ is a real Banach space with complexification $X'$ 
	 and $M \subset X$ is open, then we can take $M' \subset X'$ open,
	 and both $M'$ and $N'$ symmetric with respect to conjugation.  
	  Furthermore, if $N$ is split, then so is $N'$.
}

\biz{Proof.}{Let $(U_\alpha, \phi_\alpha :U_\alpha \rightarrow X_\alpha)$
	be a real analytic atlas on $M$ such that 
	$\phi_\alpha (U_\alpha \cap N) \subset Y_\alpha,$ a subspace of $X_\alpha$.   

	 Let $M'$ be the complexification of $M$ as in \t\tCB.  
	 After a refinement, we can complexify $(U_\alpha,\phi_\alpha)$ to
	$(U_\alpha',\psi_\alpha)$ from the same theorem.  
	 After shrinking $M'$, we may assume that every chart intersects $N$.
	 Let $Y_\alpha'$ be the complexification of $Y_\alpha$.  
	 Define
$$
	N'=\bigcup \psi_\alpha^{-1} (Y_\alpha').
$$
	 We check that this definition agrees on overlaps of $(U_\alpha')$, i.e., that
	$U_\alpha' \cap \psi_\beta^{-1}(Y_\beta')=U_\beta' \cap \psi_\alpha^{-1}(Y_\alpha')$.  
	 Choose $x \in U_\alpha' \cap \psi_\beta^{-1}(Y_\beta')$.  
	 Then $\psi_\beta (x) \in \psi_\beta (U_\alpha' \cap U_\beta') \cap Y_\beta'$.
	 Define $\phi_{\beta \alpha}= \phi_\beta \circ \phi_\alpha^{-1}$.  
	 Since $\phi_{\beta \alpha}|Y_\alpha$ is bi-\ra,
	it extends to a biholomorphic $\tilde{\phi}_{\beta \alpha}$ on a complex
	neighborhood $V\subset \psi_\alpha(U_\alpha')$ of
	$\phi_\alpha(U_\alpha \cap U_\beta) \cap Y_\alpha$.  
	 Furthermore, $\tilde{\phi}_{\beta \alpha}(V) \subset Y_\beta'$.  
	 Let $\psi_{\alpha \beta}=\psi_\alpha \circ \psi_\beta^{-1}$ (where defined) 
	for each $\alpha, \beta$.
	 Observe that by uniqueness, $\psi_{\beta \alpha}|V=\tilde{\phi}_{\beta \alpha}$.  
	 Therefore,
$$
	\psi_{\beta \alpha}|Y_\alpha':
	\psi_\alpha(U_\alpha' \cap U_\beta') \cap Y_\alpha'
	\to \psi_\beta (U_\alpha' \cap U'_\beta) \cap Y_\beta'
 $$
	is a biholomorphism. 
	 Thus $\psi_\alpha(x)=\psi_{\beta \alpha}^{-1} (\psi_\beta(x)) \in Y_\alpha'$,
	and therefore $x \in U_\beta \cap \psi_\alpha^{-1}(Y_\alpha')$.
	 The other inclusion is achieved by interchanging $\alpha$ and $\beta$.  
	 Thus, $N'$ is a complex submanifold with charts $(U'_\alpha, \psi_\alpha)$.
	 Since for every $x\in N\cap U_\alpha,$ $T_x^{1,0} N'$ is 
	the complexification of $T_x N$, $N$ is maximally real in $N'$.
	 If there is some real Banach space $X$ such that $M \subset X$
	is open, then we may use [\rS, Proposition 6.1] instead of
	\t\tCB\ to produce an $M' \subset X'$ open. 
	 We then proceed as above.  
	 We can assume that $M'$ is symmetric with respect to \cpx\ conjugation by
	intersecting it with its conjugate.  
	 If $N'$ were not symmetric with respect to conjugation,
	then its conjugate would be another possible complexification of $N$.  
	 But near $N,$ this complexifiation is unique, as we have seen above.  
	 It is easy to see that $N'$ is split if $N$ is, and in fact its
	tangent space and complement are just the complexifications of
	the tangent and complement spaces of $N$.  
	 The proof of \t\tCC\ is complete.
}

	 We note that the complexification of \ra\ \fdml\ paracompact Hausdorff \mfd{}s
	was done in [\rWB, \S\,1] by Whitney and Bruhat, who conclude their paper
	by pointing out two further articles that achieve the same.

	 Without details or proofs we remark the following.
	 If a closed split \ra\ \bsmfd\ $M$ of an open subset $\Oa$
	of a real \bspc\ $X$ admits a \ra\ \nbd\ retraction $r:\oa\to M$ from
	an open \nbd\ $\oa$ of $M$ with $M\subset\oa\subset\Oa$,
	then $r$ can be used to construct a complexification of $M$.
	 This is the case if $X$ is a separable \hspc, or, more generally,
	if $\Oa$ is a \ra\ separable \hmfd\ that admits a \ra\ Hermitian metric.
	 Indeed, then the tangent bundle $TM$ has a natural \ra\ complement
	in the restriction of $T\Oa$ to $M$, namely, the orthogonal complement
	(or metric normal bundle) $(TM)^\perp$.
	 Hence a \ra\ \nbd\ retraction $r:\oa\to M$ can be constructed in
	the usual manner, and it goes under the name of
	\kiemel{the nearest point mapping}.

	 For very special \bmfd{}s $M$ the tangent bundle $TM$ itself can serve as a
	complexification of $M$.
	 This is the case when $M$ is an affine \bmfd, i.e., $M$ admits an
	atlas of affinely related charts in the sense that the second derivatives
	of the transition functions are identically zero.
	 For, in the total space $TM$ of the tangent bundle of an affine \bmfd\ $M$
	we can define a \cpx\ structure in the same way as for $M$ a \bspc,
	which then makes $TM$ into an affine \cpx\ \bmfd\ and a complexification of $M$.

\alcim{\sD. THE PROOF OF THEOREM~\tAA.}

	 In this section we finish the proof of \t\tAA.

	 To prove (a) about the existence of \ra\ \nbd\ retractions, let
	$M'\subset\Oa'\subset X'$ be a complexification of
	$M\subset\Oa\subset X$ as in \t\tCC.
	 We may assume that both $M'$ and $\Oa'$ are symmetric \wrt\ \cpx\
	conjugation, and by \t\tBB(a) we may further assume that
	$\Oa'$ is \pscx, and $M'\subset\Oa'$ is closed in $\Oa'$.
	 \t\tBA(a) yields a real-type \holo\ retraction $r':\oa'\to M'$.
	 Setting $\oa=\oa'\cap X$ and $r=r'|\oa$, i.e, restricting to the
	real part, completes the proof \t\tAA(a) since if $x\in M\subset M'$,
	then $r(x)=x\in M$, and if $x\in\oa$, then $x=\bar x$ and
	$\overline{r(\bar x)}=r(x)$ imply that $r(x)\in M'\cap X=M$.

	 To prove (b) about the acyclicity of $\AA^E$ for a \ra\ \bvbdl\ $E\to M$,
	we can reduce (b) to its special case $M=\oa$ by pulling back $E$ by $r$
	of (a), which case in turn is a special case of \t\tBB(b).

	 To prove (c) we reduce it to its special case $M=\oa\subset X$ open by a
	\ra\ \nbd\ retraction $r:\oa\to M$ from (a).
	 To prove (c) for $M=\oa$ open, let $\oa'\subset X'$ be a symmetric
	complexification of $\oa\subset X$ with $\oa'$ \pscx\ open so small that
	the \ra\ \bvbdl\ $E\to\oa$ has for its \holo\ form the \holo\ \bvbdl\
	$E'\to\oa'$ considered by [\rS].
	 As the cocycle $(g_{UV})$ of $E'$ is \holo\ of real-type, and as
	$E'\to\oa'$ is \cts{}ly trivial, an application of \t\tBA(b) concludes
	the proof of \t\tAA(c) since the cocycle of $E$ and $E'$ over $\oa$ are
	both the same real $g_{UV}(x)$ for $x\in U\cap V\cap X$ real, and
	the trivialization $d_U(x)$ is real for $x\in U\cap X$ real, $d_U$ being
	of real-type.
	 (Note in passing that, considering, e.g., a M\"obius band over a circle, it is not
	enough to prove that the complexification $E\otimes\CC$ of a \vbdl\ $E$ is
	trivial in order to conclude that $E$ is trivial; one needs to show 
	that $E\otimes\CC$ has a trivialization of real type.)

	 Part (d) is a special case of (c) on noting that $E\to M$ is \cts{}ly trivial
	by Kuiper's theorem [\rK] stating that $\GL(\ell_2)$ is contractible.

	 Part (e) also follows from (c) since again $E\oplus(M\times Z_1)\to M$ is
	\cts{}ly trivial by [\rPC].

	 The proof of \t\tAA\ is complete.
	 
\alcim{\sE. REAL ANALYTIC DOLBEAULT ISOMORPHISM.}

	 In this section we show that the \ra\ Dolbeault isomorphism as in \S\,\sB\
	holds over certain \cpx\ \bmfd{}s.

\tetel{\t\tEA.}{The \ra\ Dolbeault isomorphism holds over a \cpx\ \bmfd\ $M$ if\/ {\rm(a)}
	or\/ {\rm(b)} below holds.
\vskip0pt
	{\rm(a)} The sheaf cohomology groups $H^q(M,\AA^E)$ for \ra\ sections
	vanish for any $q\ge1$ and any \ra\ \bvbdl\ $E\to M$.
\vskip0pt
	{\rm(b)} Our $M$ is bi-\ra\ to a closed split \ra\ \bsmfd\ $M'$ of an
	open subset\/ $\Oa$ of a real \bspc\ $X$ with an \ubs.
}

\biz{Proof.}{Part (a) follows from the formal de Rham \th\ of sheaf theory applied
	to the usual \ra\ $\delbar$-resolution of the sheaf $\OO^E$, where $E\to M$
	is any \holo\ \bvbdl.
	 Indeed, we only need to check that this resolution is locally exact, which follows
	from the local solvability of \ra\ $\delbar$-equations proved in [\rLA],
	and that the sheaf cohomology groups $H^q(M,\AA^{E^{0,p}})$ vanish for
	\ra\ $E$-valued $(0,p)$-forms for all $p\ge0$ and $q\ge1$, which in turn
	is true by our assumption.

	 Part (b) follows from (a) since $\AA^E$ is acyclic over $M$ (or $M'$)
	by \t\tAA(b), completing the proof of \t\tEA.
}

	 Note that if the condition of \t\tEA(a) holds over a \ra\ \bmfd\ $M$,
	then the \ra\ de Rham \th\ holds over $M$, i.e.,
	the sheaf cohomology group $H^q(M,\RR)$ is naturally isomorphic to
	the \ra\ de Rham group $H^q_{d,\oa}(M)$ of \ra\ real $q$-forms over $M$.

	 Grauert's famous embedding \th\ for \fdml\ \ra\ \mfd{}s (also proved for
	\cpt\ \ra\ \mfd{}s by Morrey) states that any \fdml\ \ra\ \mfd\ $M$ can be
	properly \ra{}ally embedded in \fdml\ Euclidean space, i.e., \t\tEA(b) applies
	with $\Oa=X=\RR^n$ with $n\ge1$ large enough to any \fdml\ \ra\ \mfd\ $M$,
	and thus the \ra\ Dolbeault isomorphism holds over $M$;
	a fact long known.
	 While it seems unknown whether, say, all separable \ra\ \hmfd{}s can be
	\ra{}ally embedded in \hspc, there are many abstract \ra\ \hmfd{}s that can.
	 Below are some examples.

\tetel{\p\tEB.}{{\rm(a)} The projectivization $P(\ell_2)$ of the separable \cpx\
	\hspc\ $\ell_2$ can be \ra{}ally embedded as a closed (split) \ra\ \hsmfd\
	of $\ell_2$.
\vskip0pt
	{\rm(b)} Let $M_1$ be a \cpt\ smooth \mfd,
	$M_2$ a \ra\ \mfd, and $L=C^1(M_1,M_2)$ the space of $C^1$-smooth
	maps $x:M_1\to M_2$.
	 Then $L$ can be \ra{}ally embedded as a closed split \ra\ \bsmfd\ of the
	real \bspc\ $C^1(M_1,\RR^n)$ for $n\ge1$ high enough.
}

\biz{Proof.}{To prove (a) let $x=(x_n)_{n=1}^\infty\in\ell_2\setminus\{0\}$ be a
	point $[x]$ in $P(\ell_2)$.
	 Two points $x,x'\in\ell_2$ represent the same point in $P(\ell_2)$
	if and only if there is a $c\in\CC\setminus\{0\}$ with $x'_n=cx_n$ for all
	$n\ge1$.
	 Let $y=(y_{kl})_{k,l=1}^\infty\in\ell_2$ be a point of $\ell_2$.
	 Defining $f:P(\ell_2)\to\ell_2$ by $f(x)=y$, where
$$
	y_{kl}=\frac{x_k\bar x_l}{\sum_{n=1}^\infty|x_n|^2}
 $$
 	for $k,l\ge1$, does the job.

	 This $f$ in essence is the same as the map that assigns to a point $x\not=0$
	or to a line $\CC x$ the orthogonal projector 
	$y=P_x=\frac{\langle\cdot,x\rangle x}{\|x\|^2}$ onto that line.

	 To prove (b) it is enough to take a \ra\ embedding $i:M_2\to\RR^n$ for some
	$n\ge1$ and a \ra\ retraction $r:\Oa_2\to i(M_2)$ from an open \nbd\ $\Oa_2$
	of $i(M_2)$ in $\RR^n$, and compose $i$ with the maps $x\in L$
	to get a \ra\ embedding $I:L\to C^1(M_1,\RR^n)$, $I(x)=i\circ x$.
	 The composition map $R:C^1(M_1,\Oa_2)\to I(L)$, $R(x)=r\circ x$, is a
	\ra\ retraction from an open \nbd\ $C^1(M_1,\Oa_2)$ of $I(L)$ to $I(L)$,
	thus $I(L)$ is a closed split \ra\ \bsmfd\ of $C^1(M_1,\RR^n)$.

	 The proof of \p\tEB\ is complete.
}

	  Note that in \t\tEB\ and its proof the class $C^1$ can be replaced by the
	 class $C^s$ for $s=1,2,3,\ldots$ or by the Sobolev class $W^{(s)}_2$ 
	 for $s$ positive integer large enough of maps with $s$ classical \cts\
	 or $L_2$-Sobolev derivatives.

\tetel{\t\tEC.}{{\rm(a)} The \ra\ Dolbeault isomorphism holds over any closed \cpx\ \hsmfd\
	of (any open subset of) the separable affine or projective \cpx\ \hspc.
\vskip0pt
	{\rm(b)} Any \idml\ closed \ra\ \hsmfd\ $M$ of (any open subset of)
	the separable affine or projective \hspc\ is \ra{}ally parallelizable, i.e.,
	its tangent bundle $TM$ is \ra{}ally trivial (and so are its associated
	tensor bundles, e.g., the cotangent bundle $T^*M$, the bundles $\La_p M$,
	$S_pM$ of alternating or symmetric $p$-forms on $M$).
\vskip0pt
	{\rm(c)} Let $L$ be the loop space of all loops $x:S^1\to\PP_1$ from the circle $S^1$
	to the \cpx\ projective line\/ $\PP_1$ of Sobolev class $W^{(s)}_2$, i.e.,
	with $s=1,2,\ldots$ derivatives in $L_2$, and endow $L$ with its usual structure
	of \cpx\ \hmfd.
	 Then the smooth Dolbeault group $H^{0,1}_{\delbar,\infty}(L)$, the \ra\ Dolbeault
	group $H^{0,1}_{\delbar,\oa}(L)$, and the sheaf cohomology group $H^1(L,\OO)$
	are canonically isomorphic.
	 In particular, for any $\delbar$-closed $f\in C^\infty_{0,1}(L)$ there are
	a $\delbar$-closed form $g\in C^\oa_{0,1}(L)$ and a \fn\ $h\in C^\infty(L)$ with
	$f=g+\delbar h$ on $L$.
}

\biz{Proof.}{(a) By \p\tEB(a) we may embed our \hmfd\ in the affine space $\ell_2$, and
	apply \t\tEA(b) to conclude the proof.

	Part (b) follows from \t\tAA(d).

	Part (c) follows from Ning Zhang's Dolbeault isomorphism of $H^{0,1}_{\delbar,\infty}(L)$
	and $H^1(L,\OO)$ together with our \ra\ Dolbeault isomorphism \t\tEA(b) noting that
	$L$ \ra{}ally embeds into \hspc\ as a closed \ra\ \hsmfd\ by \p\tEB(b).
}

\alcim{\sF. SMOOTH FORMS REAL ANALYTIC OUTSIDE A BALL.}

	 In this section we give a partial generalization of
	Ligocka's form of a well-known theorem of Ehrenpreis on
	the $\delbar$-equation with \cpt\ support.

	 Ehrenpreis [\rE, \rH] showed that if $f$ is
	a $\delbar$-closed $(p,q)$-form on $\CC^n$ that is zero outside a
	ball, then the $\delbar$-equation $\delbar u=f$ has a solution on $\CC^n$
	(with \cpt\ support if $q=1$).
	 Ligocka [\rLg] generalized the above \th\ to \bspc{}s.
	 In what follows we replace the support condition of vanishing outside a ball
	by a condition on real analyticity outside a ball.
	 See \t\tFA.

\tetel{\t\tFA.}{Let $X$ be a \cpx\ \bspc\ with an \ubs\ and $C^\infty$-smooth cutoff
	\fns\ (e.g., $X=\ell_2$ or $c_0$), $f\in C^\infty_{p,q}(X,Z)$, where
	$p\ge0$, $q\ge1$ and $Z$ is any \bspc.
	 If $\delbar f=0$ on $X$ and $f(x)$ is \ra\ for\/ $\|x\|>1$, then there is
	a $u\in C^\infty_{p,q-1}(X,Z)$ with $\delbar u=f$ on $X$.
}

	 Very little is known about the $\delbar$-equation, say, on a \hspc.
	 It seems unknown, e.g., whether a $C^\infty$-smooth $\delbar$-equation is
	locally solvable on $X=\ell_2$.
	 On the positive side, the only results that seem to be known 
	(especially for $q\ge2$) are the \th{}s of Lempert in [\rLA]
	that derive from the projective space and impose a global growth rate, 
	the \th\ of Ehrenpreis and Ligocka, 
	that \ra\ $\delbar$-equations are globally solvable on $X$,
	and the intermediate result between the last two, expressed in \t\tFA.

	 As we shall use the \ra\ Dolbeault isomorphism as in \t\tBB(c),
	we study the sheaf cohomology groups $H^q(X\setminus B,\OO^Z)$ first,
	where $B\subset X$ is a closed ball, or in fact a certain closed convex set.

\tetel{\p\tFB.}{Let $Y$ be a \cpx\ \bspc\ with an \ubs, $\Da\subset\CC^n$ a \cpt\
	polydisc of nonempty interior, $n\ge3$, $F\subset\OO(Y)$ nonempty,
	$K=\{y\in Y:|f(y)|\le1\text{\ for all\ }f\in F\}$, $L=\Da\times K$, $X=\CC^n\times Y$,
	$\Oa=X\setminus L$, and $Z$ any \bspc.
	 Then $H^q(\Oa,\OO^Z)$ vanishes for all\/ $1\le q\le n-2$.
}

\biz{Proof.}{Assume as we may that $\Da$ is the unit polydisc
	$\Da=\{z\in\CC^n:|z_i|\le1\text{\ for all\ }i=1,\dots,n\}$.
	 We use Laurent series expansions \wrt\ the $z_i$, $i=1,\ldots,n$,
	and the Leray covering method of Frenkel's lemma as in
	[\rC, \rF, \rG, \rJE, \rPB].
	 Let $X\ni x=(z,y)\in\CC^n\times Y$,
	$U_i=\{x\in X:|z_i|>1\}$ for $i=1,\ldots,n$,
	$U_f=\{x\in X:|f(y)|>1\}$ for $f\in F$, and
	$\UU=\{U_i,U_f:i=1,\ldots,n; f\in F\}$.
	 Then each $U\in\UU$ is \pscx\ open in $X$ and is a product set of the form
	$U=A_1\times\ldots\times A_n\times B$, where $A_i\subset\CC$ for $i=1,\ldots,n$,
	and $B\subset Y$, and thus so is the intersection of any finitely many members
	of $\UU$, i.e., the covering $\UU$ of $\Oa$ is a Leray covering for the sheaf
	$\OO^Z$ over $\Oa$ by the vanishing theorem of Lempert [\rLC], hence
	$H^q(\Oa,\OO^Z)\cong H^q(\UU,\OO^Z)$ for all $q\ge1$, and thus we only have
	to show that $H^q(\UU,\OO^Z)=0$ for $1\le q\le n-2$.

	 If $f\in\OO(V,Z)$ is a \holo\ \fn\ on an open set $V\subset X$ of the
	form $V=A_1\times\ldots\times A_n\times B$, where $A_i\subset\CC$ is an
	annulus for $i=1,\ldots,n$, then $f$ can be expanded in Laurent series
	with positive and negative Laurent projections $P^+_if$ and $P^-_if$ given by
$$\eqalign{
	f(x)&=
	  \sum_{k=-\infty}^\infty f_{i,k}(z_1,\ldots,z_{i-1},z_{i+1},\ldots,z_n,y) z_i^k\cr
        (P_i^+f)(x)&=
          \sum_{k=0}^\infty f_{i,k}(z_1,\ldots,z_{i-1},z_{i+1},\ldots,z_n,y) z_i^k,\cr
        (P_i^-f)(x)&=
          \sum_{k=-\infty}^{-1} f_{i,k}(z_1,\ldots,z_{i-1},z_{i+1},\ldots,z_n,y) z_i^k\cr
}
 $$
 	for $x\in V$ for any $i=1,\ldots,n$, and each coefficient $f_{i,k}\in\OO(V)$
	is a \holo\ \fn\ at least on $V$.

	 If $i=1,\ldots,n$, $V=\bigcap_{j=0}^q V_j$, $V_j\in\UU$, and no $V_j$ equals $U_i$,
	then $P_i^-f=0$ since $f=P_i^+f$ is an entire \fn\ \wrt\ $z_i$.
	 If some $V_j$ equals $U_i$, then $P^+_if$ is entire in $z_i$, and thus
	$P^+_if\in\OO(\bigcap_{j\not=i}V_j,Z)$.
	 Looking, for $i=1,\ldots,n$, at the chain map $f\mapsto Qf=g$,
	where $f=(f_{V_0\ldots V_q})\in C^q(\UU,\OO^Z)$, and
	$g=(g_{V_0\ldots V_{q-1}})\in C^{q-1}(\UU,\OO^Z)$ is defined by
	$g_{V_0\ldots V_{q-1}}=P^+_if_{U_iV_0\ldots V_{q-1}}$,
	it can be seen as in [\rG] that 
	$P^-_if=f-\da Qf-Q\da f$, and so
	any cocycle $f$ is cohomologous to the
	cocycle $P_i^-f$ in $H^q(\UU,\OO^Z)$ for all $i=1,\ldots,n$.
	 Hence $f$ and $h=P_1^-\ldots P_n^-f$ are also cohomologous in $H^q(\UU,\OO^Z)$.
	 If $h_{V_0\ldots V_q}\not=0$, then each of $U_1,\ldots,U_n$ must appear among
	$V_0,\ldots,V_q$.
	 In particular, $n\le q+1$.
	 Thus $h$ vanishes for $1\le q\le n-2$, and so do $H^q(\UU,\OO^Z)$ and
	$H^q(\Oa,\OO^Z)$, completing the proof of \p\tFB.
}

	 Note that the Hahn--Banach \th\ tells us that any closed convex set $K\subset Y$
	can be represented in the form required in \p\tFB\ with $F\subset\OO(Y)$,
	where $f(y)=e^{\eta_f(y)}$, and $\eta_f\in Y^*$ is a \cts\ linear functional
	for all $f\in F$.
	 In particular, $K$ can be any closed ball in $Y$.

\tetel{\p\tFC.}{{\rm(a)} Let $X$ be a para\cpt\ Hausdorff space, $S\to X$ a sheaf
	of Abelian groups, $\oa\subset X$ open, $L\subset X$ closed with $L\subset\oa$.
	 If the sheaf cohomology groups $H^q(X,S)$ and $H^q(\oa,S)$ vanish for all $q\ge1$,
	then the groups $H^q(\oa\setminus L,S)$ and $H^q(X\setminus L,S)$ are isomorphic
	for all $q\ge1$.
\vskip0pt
	{\rm(b)} If $X,Z,L$ are as in \p\tFB, $\oa\subset X$ \pscx\ open with $L\subset\oa$,
	then $H^q(\oa\setminus L,\OO^Z)$ vanishes for all $1\le q\le n-2$.
}

\biz{Proof.}{Part (a) follows directly from the Mayer--Vietoris long exact sequence
$$
\eqalign{
	0\to H^0(X,S)\to H^0(\oa,S)\oplus H^0(X\setminus L,S)\to H^0(\oa\setminus L,S)\to\cr
	\phantom{0\to }H^1(X,S)\to H^1(\oa,S)\oplus H^1(X\setminus L,S)\to H^1(\oa\setminus L,S)\to\cr
	\ldots\phantom{a nice long gap}\ldots\phantom{a nice long gap}
	\ldots\phantom{a nice long gap}\ldots\cr
	\phantom{0\to }H^q(X,S)\to H^q(\oa,S)\oplus H^q(X\setminus L,S)\to H^q(\oa\setminus L,S)\to\cr
	\ldots\phantom{a nice long gap}\ldots\phantom{a nice long gap}
	\ldots\phantom{a nice long gap}\ldots\cr
}
 $$
 	of sheaf cohomology applied to the sheaf $S$ and the couple $\oa$ and $X\setminus L$
	that forms an open covering of $X$,
	while (b) follows immediately from (a) and \p\tFB, completing the proof of \p\tFC.
}

	 Using \p\tFC(b) one can also formulate a \th\ analogous to \t\tFA\ 
	that states that certain smooth
	$\delbar$-equations on a \pscx\ open set can be solved if the given form
	is \ra\ in a suitable inner collar of the boundary of the domain.

\biz{Proof of \t\tFA.}{Fix $q\ge1$ and choose $n\ge3$ with $q\le n-2$, a decomposition
 	$X=\CC^n\times Y$ of \bspc{}s, and a \cpt\ polydisc $\Da\subset\CC^n$ and a closed
	\bdd\ ball $K\subset Y$ so that our form $f$ is \ra\ on $\Oa=X\setminus L$,
	where $L=\Da\times K$.
	 Choose a \cpt\ polydisc $\Da'\subset\CC^n$ that contains $\Da$ in its interior,
	and a closed \bdd\ ball $K'\subset Y$ that contains $K$ in its interior,
	and let $L'=\Da'\times K'$, and $\Oa'=X\setminus L'$.
	 As $H^q(\Oa,\OO^Z)=0$ by \p\tFB\ we see by the \ra\ Dolbeault isomorphism
	\t\tBB(c) that $H^{p,q}_{\delbar,\oa}(\Oa,Z)=0$, so there is a \ra\ $(p,q)$-form
	$v\in\AA_{p,q}(\Oa,Z)$ with $\delbar v=f$ on $\Oa$.
	 Let $\chi\in C^\infty(X)$ be a smooth cutoff \fn\ that equals $1$ on a \nbd\ of
	the closure of $\Oa'$ and $0$ on a \nbd\ of $L$.
	 Look at $g=f-\delbar(\chi v)\in C^\infty_{p,q}(X,Z)$.
	 Then $g$ is $\delbar$-closed on $X$ and $g=0$ on $\Oa'$.
	 The \th\ of Ehrenpreis and Ligocka yields a $u\in C^\infty_{p,q-1}(X,Z)$
	with $\delbar u=g$, hence $f=\delbar(u+\chi v)$ is $\delbar$-exact on $X$
	as claimed.
}

	 Note that in the Ehrenpreis--Ligocka \th\ we do not need $g$ to be 
	$C^\infty$-smooth but, say, $C^1$-smooth, and hence in the proof above 
	we do not need the cutoff \fn\ $\chi$ to be $C^\infty$-smooth, 
	but, say, $C^2$-smooth; hence
	$X$ could be $X=\ell_p$, $2\le p<\infty$ in \t\tFA.

\alcim{\sG. REMARKS ON REAL ANALYTIC BANACH MANIFOLDS.}

	 In this section we make some geometrical remarks on \ra\ \bmfd{}s.

\tetel{\t\tGA.}{Let $X$ be a real \bspc\ with an \ubs, $\Oa\subset X$ open,
	and $M\subset\Oa$ a closed split \ra\ \bsmfd\ of\/ $\Oa$.
	 Then $M$ is bi-\ra\ to a closed split \ra\ \bsmfd\ $N$ of a \bspc\ $Y$
	with an \ubs.
	 Moreover, if $X$ is a separable $L_p$-space, $1<p<\infty$, then
	$Y$ can also be taken to be a separable $L_p$-space for the same $p$.
}

\biz{Proof.}{Let $\Oa'\subset X'$ be a complexification of
	$\Oa\subset X$ as in \t\tBB(b).
	 An application of Zerhusen's embedding theorem [\rZn]
	to $\Oa'\subset X'$ completes the proof of \t\tGA.
}

	 We next show that covering spaces of embeddable \ra\ \bmfd{}s are
	embeddable \ra\ \bmfd{}s.

\tetel{\t\tGB.}{Let $X,\Oa,M$ be as in \t\tGA, and $\tilde M\to M$ a covering space
	of $M$ with countably many sheets.
	 Then $\tilde M$ is bi-\ra\ to a closed split \ra\ \bsmfd\ of a \bspc\
	with an \ubs.
}

\biz{Proof.}{We follow Lempert's application of the \holo\ version of \t\tAA(d) to
	\holo\ covering spaces in [\rLE, Thm.\,3.7].
	 Let $N\subset\NN$ be \st\ the bundle $\tilde M\to M$ admits trivializations
	with cocycle $(g_{UV})\in Z^1(\UU,\Aut(N))$, where $\UU$ is an open covering
	of $M$, and $\Aut(N)$ is the group of all permutations of the discrete set $N$.
	 Let $\ell_2$ be the \hspc\ with orthonormal basis $e_n$, $n\in\NN$, and
	for any $p\in\Aut(N)$ define the permutation operator $\dot p\in\GL(\ell_2)$
	by $\dot pe_n=e_{p(n)}$ for $n\in N$ and $\dot pe_n=e_n$ for $n\in\NN\setminus N$.
	 Clearly, $\dot p$ is an isometry of $\ell_2$.
	 Look at the `flat' \ra\ \hvbdl\ $E$ defined by the \ra\ cocycle
	$(\dot g_{UV})\in Z^1(\UU,\AA^{\GL(\ell_2)})$, which is just the permutation
	representation of the original cocycle $(g_{UV})$.
	 There is a \ra\ embedding $i:\tilde M\to E$ into the total space $E$ of this
	\ra\ \hvbdl.
	 As there is a \ra\ isomorphism $j:E\to M\times\ell_2$ by \t\tAA(d), and as
	$j(i(\tilde M))$ is a closed split \ra\ \bsmfd\ of $M\times\ell_2$ (as it is
	easily seen in local trivializations), an application of \t\tGA\ to $M\times\ell_2$
	completes the proof of \t\tGB.
}

\tetel{\t\tGC.}{Let $K$ be a \ra\ \cpt\ \mfd, and $M,\Oa,X$ as in \t\tGA.
	 Then any \cts\ map $f:K\to M$ is homotopic to a \ra\ map $g:K\to M$.
}

\biz{Proof.}{Let $r:\oa\to M$ be a \ra\ retraction as in \t\tAA(a) of an open set
	$\oa$ with $M\subset\oa\subset\Oa\subset X$.
	 The Malgrange approximation \th\ or Grauert's \ra\ embedding of $K$ into
	$\RR^n$ as $K'\subset\RR^n$ for $n\ge1$ high enough followed by an application
	of the Weierstrass approximation \th\ for polynomials $K'\to X$ gives us a
	\ra\ $g':K\to X$ \st\ the linear homotopy $(1-t)f(x)+tg'(x)$ takes its values
	in $\oa$ for all $x\in K$ and $t\in[0,1]$.
	 Thus the homotopy $h:K\times[0,1]\to M$ defined by $h(x,t)=r((1-t)f(x)+tg'(x))$
	makes sense and completes the proof of \t\tGC\ on letting $g(x)=h(x,1)=r(g'(x))$
	for $x\in K$.
}

	 In particular, the homotopy classes of spheroids in the homotopy groups
	$\pi_n(M)$, $n\ge1$, can be represented by \ra\ maps.

	 In conclusion we remark that the method of complexification as treated in
	this paper can be used to obtain further information about \ra\ \bmfd{}s, and
	\ra\ \fns\ and sections over them, e.g., in the manner of [\rDPV].

\vskip0.30truein
\centerline{\scVIII References}
\vskip0.20truein
\baselineskip=11pt
\parskip=3pt
\frenchspacing
{\rmVIII

	[\rC] Cartan,~H.,
	{\itVIII
	Sur le premier probl\`eme de Cousin},
	C. R. Acad. Sci. Paris, {\bfVIII 207} (1938), 558--560.

	[\rD] Dineen, S.,
	{\itVIII 
	Complex analysis in infinite dimensional spaces},
	Springer-Verlag, 
	London,
	(1999).

	[\rDPV] \vonal, Patyi, I., Venkova, M.,
	{\itVIII 
	Inverses depending holomorphically on a parameter in a Banach space},
	J. Funct. Anal., 
	{\bfVIII 237} (2006), no. 1, 338--349.

	[\rE] Ehrenpreis,~L.,
	{\itVIII
	A new proof and an extension of Hartogs' theorem},
	Bull. Amer. Math. Soc.,
	{\bfVIII 67} (1961), 507--509.

	[\rF] Frenkel,~J.,
	{\itVIII
	Cohomologie non ab\'elienne et espaces fibr\'es},
	Bull. Soc. Math. France, {\bfVIII 85} (1957), 135--220.

        [\rG] Gunning,~R.C.,
        {\itVIII
        Introduction to holomorphic functions of several variables},
        Vol. III,
        Wadsworth \& Brooks/Cole, Belmont, California,
        (1990).

	[\rJE] Honda,~T., Miyagi,~M., Nishihara,~M., Ohgai,~S., Yoshida,~M.,
	{\itVIII 
	The Frenkel's lemma in Banach spaces and its applications},
	Far East J. Math. Sci., 
	{\bfVIII 14} (1) (2004), 69--93.

        [\rH] H\"ormander,~L.,
        {\itVIII  
        An Introduction to Complex Analysis in Several Variables},
        3rd Ed., North-Holland, Amsterdam,
        (1990).

	[\rK] Kuiper,~N.H.,
	{\itVIII
	The homotopy type of the unitary group of Hilbert space},
	Topology,
	{\bfVIII 3}
	(1965),
	19--30.

	[\rLA] Lempert,~L.,
	{\itVIII
	The Dolbeault complex in infinite dimensions~I},
	J. Amer. Math. Soc., {\bfVIII 11} (1998), 485--520.

        [\rLC] \vonal,
        {\itVIII
        The Dolbeault complex in infinite dimensions~III},
        Invent.{} Math., {\bfVIII 142} (2000), 579--603.

        [\rLD] \vonal,
        {\itVIII
        Vanishing cohomology for holomorphic vector bundles in a Banach setting},
        Asian J. Math.,
        {\bfVIII 8}
        (2004),
        65--85.

	[\rLE] \vonal,
	{\itVIII
	Analytic continuation in mapping spaces},
	manuscript.

	[\rLP] \vonal, Patyi,~I.,
	{\itVIII
	Analytic sheaves in Banach spaces},
	Ann. Sci. \'Ecole Norm. Sup., 
	S\'er. 4,
	{\bfVIII 40}
	(2007),
	453--486.
	
	[\rLg] Ligocka,~E.
	{\itVIII 
	Levi forms, differential forms of type (0,1) and pseudoconvexity in Banach spaces},
	Ann. Polon. Math.,
	{\bfVIII 33} (1976), no. 1-2, 63--69. 

	[\rM] Mujica, J.,
	{\itVIII
	 Complex analysis in Banach spaces},
	North--Holland, Amsterdam,
	(1986).

        [\rPA] Patyi,~I.,
        {\itVIII
        On the Oka principle in a Banach space I},
        Math. Ann.,
        {\bfVIII 326}
        (2003),
        417--441.

	[\rPB] \vonal,
	{\itVIII
	Cohomological characterization of pseudoconvexity in a Banach space},
	Math.Z., {\bfVIII 245} (2003), 371--386. 

        [\rPC] \vonal,
        {\itVIII
        On holomorphic Banach vector bundles over Banach spaces},
        manuscript.

	[\rS] Simon,~S.B.,
	{\itVIII 
	A Dolbeault isomorphism theorem in infinite dimensions},
	Trans. Amer. Math. Soc., to appear.

	[\rWB] Whitney,~H., Bruhat,~F.,
	{\itVIII Quelques propri\'et\'es fondamentales des ensembles anal\-ytiques-r\'eels},
	Comment. Math. Helv.,
	{\bfVIII 33}, (1959), 132--160.

	[\rZn] Zerhusen,~A.B.,
	{\itVIII
	Embeddings of pseudoconvex domains in certain Banach spaces},
	Math. Ann.,
	{\bfVIII 336}
	(2006), no. 2,
	269--280.

	[\rZg] Zhang,~N.,
	{\itVIII
 	The Picard group of a loop space},
	manuscript, arxiv:math.CV/0602667.

}
\vskip0.20truein
\centerline{\vastag*~***~*}
\vskip0.15truein
{\scVIII
        Imre Patyi,
        Department of Mathematics and Statistics,
        Georgia State University,
        Atlanta, GA 30303-3083, USA,
        {\ttVIII ipatyi\@gsu.edu}

	Scott B. Simon,
	Department of Mathematics,
	SUNY Stony Brook,
	Stony Brook, NY 11794-3651, USA,
	{\ttVIII sbsimon\@math.sunysb.edu}
}
\bye